\documentclass[11pt]{amsart}
\usepackage{amsfonts,latexsym,rawfonts,amsmath,amssymb,amsthm,graphicx}

\numberwithin{equation}{section}
\newtheorem{theorem}{Theorem}[section]
\newtheorem{lem}[theorem]{Lemma}
\newtheorem{thm}[theorem]{Theorem}

\newtheorem{cor}[theorem]{Corollary}

\newtheorem{rem}[theorem]{Remark}

\def\s{\,\,\,\,}

\def\dis{\displaystyle}
\def\dint{\displaystyle{\int}}

\def\endproof{$\hfill\Box$\\}

\def\R{\mathbb{R}}

\def\C{\mathbb{C}}

\title{\bf An example of $\alpha$-harmonic map sequence for which energy identity is not true}
\author[Y. Li \& Y. Wang]{Yuxiang Li\\
{\small\it Department of Mathematical Sciences},\\
{\small\it Tsinghua University,}\\
{\small\it Beijing 100084, P.R.China.}\\
{\small\it Email: yxli@math.tsinghua.edu.cn.}\\ \\
Youde Wang\\
{\small\it Academy of Mathematics and Systems Sciences,}\\
{\small\it Chinese Academy of Sciences,}\\
{\small\it Beijing 100080,  P.R. China.}\\
{\small\it Email: wyd@math.ac.cn}}

\date{}
\begin{document}
\maketitle

\begin{abstract}
We construct a closed Riemannian manifold $(N,h)$ and a sequence of $\alpha$-harmonic
maps from $S^2$ into $N$ with uniformly bounded energy such that the energy identity for this
sequence is not true.

\end{abstract}

\section{\bf Introduction}
Let $(\Sigma,g)$ be a Riemann surface and $(N,h)$ be an $n$-dimensional smooth compact Riemannian
manifold which is embedded in $\mathbb{R}^K$. Usually, we denote the space of Sobolev maps from $\Sigma$
into $N$ by $W^{k,p}(\Sigma, N)$, which is defined by
$$W^{k,p}(\Sigma, N)=\{u\in W^{k,p}(\Sigma, \mathbb{R}^K): u(x)\in N\,\, \text{for a.e.}\,\,x\in\Sigma\}.$$
For $u\in W^{1,2}(\Sigma,N)$, we define locally the energy
density $e(u)$ of $u$ at $x\in \Sigma$ by
$$e(u)(x)=|\nabla_g u|^2=g^{ij}(x)h_{\alpha\beta}(u(x))
\frac{\partial u^\alpha}{\partial x^i}\frac{\partial
u^\beta}{\partial x^j}.$$
The energy of $u$ on $\Sigma$, denoted by $E(u)$ or $E(u, \Sigma)$, is defined by
$$E(u)=\frac{1}{2}\dint_\Sigma e(u)dV_g,$$
and the critical points of $E$ are called harmonic maps. We know
that a harmonic map $u$ satisfies the following equation:
$$\tau(u)=\Delta u+A(u)(\nabla u,\nabla u)=0,$$
where $A$ is the second fundamental form of $N$ in $\mathbb{R}^K$.
Harmonic maps are related very closely to minimal surface.
It is well-known that a harmonic map from $S^2$ into $N$ must be a branched minimal
immersion in $N$.

Unfortunately, $E$ does not satisfy the Palais-Smale condition. From the viewpoint of
calculus of variation, it is difficult to show the existence of harmonic maps from a surface.
In order to obtain harmonic maps, in \cite{Sacks-Uhlenbeck1} Sacks and Uhlenbeck introduced
the so called $\alpha$-energy $E_\alpha$ instead of $L^2$ energy $E$ as the following
$$E_\alpha(u)=\frac{1}{2}\int_\Sigma\{(1+|\nabla u|^2)^\alpha-1\} dV_g,$$
where we always assume that $\alpha>1$. It is well-known that this $\alpha$-energy functional
$E_\alpha$ satisfies the Palais-Smale condition. The critical points of $E_\alpha$ in $W^{1,2\alpha}(\Sigma,N)$,
which are called as the $\alpha$-harmonic maps, satisfy the following equation:
\begin{equation}\label{alpha}
\Delta_gu_\alpha+(\alpha-1)\frac{\nabla_g|\nabla_g u_\alpha|^2
\nabla_gu_\alpha}{1+|\nabla_gu_\alpha|^2} +
A(u_\alpha)(du_\alpha,du_\alpha)=0.
\end{equation}
The strategy of Sacks and Uhlenbeck is to employ a sequence of $\alpha$-harmonic maps
to approximate a harmonic map as $\alpha$ tends to 1. Hence, to show the existence of
harmonic maps we need to study the convergence behavior of a sequence of
$\alpha$-harmonic maps $u_\alpha$ with $E_\alpha(u_\alpha) <C$ from a compact surface
$(\Sigma, g)$ into a compact Riemannian manifold $(N,h)$ without boundary. Generally,
such a sequence converges weakly to a harmonic map in $W^{1,2}(\Sigma, N)$ and strongly
in $C^\infty$ away from a finite set of points in $\Sigma$.

Concretely, let $\{u_{\alpha_k}\}$ be a sequence of $\alpha$-harmonic maps
from $\Sigma$ into $N$ with uniformly bounded $\alpha$-energy, i.e.,
$E_{\alpha_k}(u_{\alpha_k})<\Lambda<\infty$. We assume that the sequence does not converge smoothly on $\Sigma$.
By the theory of Sacks-Uhlenbeck, there exists a subsequence of $\{u_{\alpha_k}\}$, still denoted by
$\{u_{\alpha_k}\}$, and a finite set $\mathcal{S}\subset\Sigma$ such that the subsequence
converges to a harmonic map $u_0$ in $C^\infty_{loc}(\Sigma\setminus \mathcal{S})$. We know that, at each point
$p_i\in\mathcal{S}$, the energy of the subsequence concentrates and the blow-up phenomena occur.
Moreover, there exist point sequences $\{x_{i_k}^l\}$ in $\Sigma$ with $\lim\limits_{k\rightarrow+\infty}x_{i_k}^l= p_i$ and scaling
constant number sequences $\{\lambda_{i_k}^l\}$ with $\lim\limits_{k\rightarrow+\infty}\lambda_{i_k}^l\rightarrow 0$,
$l=1,\cdots, n_0$, such that
$$u_{\alpha_k}(x_{i_k}^l+\lambda_{i_k}^lx)\rightarrow v^l\s\s \text{in}\s\s C^j_{loc}(\mathbb{R}^2 \setminus \mathcal{A}^i),$$
where all $v^i$ are non-trivial harmonic maps from $S^2$ into $N$, and $\mathcal{A}^i\subset\mathbb{R}^2$ is a finite set.

In order to explore and describe the asymptotic behavior of $\{u_{\alpha_k}\}$ at each blow-up point,
the following two problems were raised naturally. One is whether or not the energy identity holds true, i.e.,
$$\lim_{\alpha_k\rightarrow 1}E_{\alpha_k}(u_{\alpha_k}, B^\Sigma_{r_0}(p_i))= E(u_0,B^\Sigma_{r_0}(p_i))+ \sum_{l=1}^{n_0}E(v^l).$$
Here, $B^\Sigma_{r_0}(p_i)$ is a geodesic ball in $\Sigma$ which contains only one blow-up point $p_i$. The other is whether or
not the necks connecting bubbles are some geodesics of finite length?

One paid attention to the two problems for long time and has made considerably
great progress. Now, let us recall some main results on the problems. Chen and Tian
in \cite{Chen-Tian} considered a special sequence $\{u_{\alpha_k}\}$ with uniformly
bounded $\alpha$-energy, for which every  $u_{\alpha_k}$ is a minimizing
$\alpha_k$-harmonic map and all maps $u_{\alpha_k}$ belong to a fixed homotopy class.
They studied the convergence behavior of such a special sequence and provided a proof
on the above energy identity. Later, for the same sequence, Li and Wang in \cite{Li-Wang2}
gave another constructing proof on the energy identity, which is completely different
from that given in \cite{Chen-Tian}.

One also considered the energy identity for a minimax sequence of
$\alpha$-harmonic maps. Suppose that $A$ is a parameter manifold. Let $h_0:\Sigma\times A
\rightarrow N$ be a continuous map, and $H$ be such a set of
continuous maps $h: \Sigma\times A \rightarrow N$ that every $h\in H$ is
homotopic to $h_0$ and satisfies $h(t)\in W^{1, 2\alpha}(\Sigma, N)$ for
any fixed $t\in A$. Set
\begin{equation*}\label{jost}
\beta_\alpha(H)=\inf_{h\in H}\sup_{t\in A}E_\alpha(h(\cdot,t)).
\end{equation*}
One has known that there is at least a sequence $\{u_{\alpha_k}\}$,
each $u_{\alpha_k}$ of which attains $\beta_{\alpha_k}(H)$,
satisfies the energy identity as $\alpha_k\rightarrow 1$.
For more details, we refer to \cite{J, Lamm}.

On the other hand, it should be pointed out that one has also established some effective methods to prove
successfully the energy identity and give the detailed description of the connecting necks
for the heat flow of harmonic maps from a Riemann surface, or more generally, a
sequence of maps from a Riemann surface with tension fields $\tau$ bounded in the
sense of $L^2$ \cite{Ding, Ding-Tian, Qing,Qing-Tian}.

Recently, Li and Wang in \cite{Li-Wang1} studied the above problems on the sequences
of $\alpha$-harmonic maps and obtained some results which can be summarized as follows.
If the energy concentration phenomena appears for $\{u_{\alpha_k}\}$,
a weak energy identity was proved and a direct convergence relation between
the blow-up radius and the parameter $\alpha$, which ensures the energy identity and no-neck property,
was discovered. Li and Wang also showed that the necks converge to some geodesics and gave a length
formula for the neck in the case only one bubble appears.

Motivated by an example given by Duzaar and Kuwert in \cite{Duzaar-Kuwert}, Li and Wang \cite{Li-Wang1}
also constructed an $\alpha$-harmonic map sequence with uniformly bounded energy, for which the blow-up
phenomenon occurs and there exists at least a neck (geodesic) of infinite length. This answers negatively
the second problem on $\alpha$-harmonic map sequence.

Although some mathematicians think that the energy identity for the sequence of $\alpha$-harmonic maps should also
be true, up to now it has been unclear in general whether the energy identity
for an $\alpha$-harmonic map sequence with bounded energy holds true or
not. In this short paper, we will modify the construction in \cite{Li-Wang1}
to show that the energy identity is also not true.\\

On the other hand, a natural problem is whether the set of the values of energy for harmonic spheres in any given Riemannian
manifold $(N, h)$ is discrete or not, since the bubbles produced in the convergence of a sequence of $\alpha$-harmonic maps
from $(\Sigma, g)$ are always harmonic spheres.

We denote the set by
$$\mathcal{E}(N,h)=\{E(u):u \mbox{ is a harmonic map from $S^2$ into $(N, h)$}\}.$$
It is well known that, if $(N,h)$ is the standard sphere $S^2$, we have
$$\mathcal{E}(N,h)=\{4k\pi: k=0, 1, \cdots, n, \cdots\}.$$
We also know from \cite{Va} due to Valli that, if $(N, h)$ is the unitary group $U(n)$ with the standard metric, then the energy of harmonic maps
$S^2\rightarrow U(n)$ can take as values only integral multiples of $8\pi$. Some other energy gap phenomenon on unitons were discussed in
\cite{Anand, Dong, Uhlenbeck}. Some mathematicians conjectured that $\mathcal{E}(N,h)$ is a discrete set.
Here, we will also give a counterexample to show that $\mathcal{E}(N,h)$ is not discrete.

\section{\bf $\alpha$-harmonic maps}
Later, we will discuss the convergence behavior of some $\alpha$-harmonic map sequences with
uniformly bounded $\alpha$-energy or $L^2$ energy. In fact, by discussing the convergence of
$\alpha$-harmonic map sequences, Sacks and Uhlenbeck developed an existence theory on minimal
surfaces in \cite{Sacks-Uhlenbeck1, Sacks-Uhlenbeck2}. In particular, they established in \cite{Sacks-Uhlenbeck1}
the well-known $\epsilon$-regularity theorem on $\alpha$-harmonic maps and removal singularity theorem
on harmonic maps, which will be used repeatedly in the present paper.
\begin{thm}\label{epsilon}
Let $D=D_1(0)=\{z: |z|<1\}\subset\mathbb{C}$ be a disk with radius $1$ and $N$ be a Riemannian manifold.
Assume that $u: D\rightarrow N$ satisfies equation \eqref{alpha}. Then, there exists $\epsilon_0>0$
and $\alpha_0>1$ such that, if $E(u,D)<\epsilon_0$ and
$1\leq\alpha\leq\alpha_0$,  then we have
$$\|\nabla^{k} u\|_{L^\infty(D_\frac{1}{2})}\leq C(k)E(u,D).$$
\end{thm}

\begin{thm}\label{removable}
Assume that $u: D\setminus\{0\}\rightarrow N$ is a
harmonic map with $E(u)<+\infty$. Then, $u$ is a harmonic map
from $D$ into $N$.
\end{thm}
The above theorem tells us that, if $u$ is a harmonic from $\C\setminus\{p_i\in\C: i=1,2,\cdots,l<\infty\}$ into $N$ with
$E(u)<+\infty$, then $u$ can be viewed as a harmonic map from $S^2$ into $N$.

Now, we can state more precisely the energy concentration of $\{u_{\alpha_k}\}$.
Let $B^\Sigma_t(x)$ denote the geodesic ball of $\Sigma$ which is centered at $x$ and of geodesic radius $t$.
By Theorem 2.1, we know that the finite singular set of $\{u_{\alpha_k}\}$ can be defined precisely by
$$\mathcal{S}=\left\{x\in\Sigma:\lim_{t\rightarrow 0}\varliminf_{k \rightarrow +\infty}\int_{B^\Sigma_t(x)}|\nabla u_{\alpha_k}|^2\geq\frac{\epsilon_0}{2}\right\}.$$
For any $\tilde{x}_0\notin \mathcal{S}$, there exists $\delta>0$ such that
$E(u_{\alpha_k},B^\Sigma_\delta(\tilde{x}_0))<\epsilon_0$. Applying Theorem \ref{epsilon}, we know that $\{u_{\alpha_k}\}$
converges smoothly on any $\Omega\subset\subset\Sigma\setminus\mathcal{S}$. The limit map is a harmonic map from
$\Sigma\setminus\mathcal{S}$ into $N$. Theorem \ref{removable} tells us that the singular points of the limit map
can be removed, i.e., it is a harmonic map from $\Sigma$ into $N$.

If $x_0\in\mathcal{S}$, it is easy to check that, for any $t$, there holds
$$\|\nabla u_{\alpha_k}\|_{C^0(B^\Sigma_t(x_0))}\rightarrow+\infty.$$
Choose $x_{\alpha_k}\in B^\Sigma_{\delta}(x_0)$ such that $$|\nabla u_{\alpha_k}(x_{\alpha_k})|=\max_{B^\Sigma_{\delta}(x_0)}|\nabla
u_{\alpha_k}|$$ and let
$$\lambda_{\alpha_k}=\frac{1}{\max_{B^\Sigma_{\delta}(x_0)}|\nabla u_{\alpha_k}|}.$$
It is easy to see that $x_{\alpha_k}\rightarrow x_0$ as $k\rightarrow\infty$.
Then, in an isothermal coordinate system around $x_0$, we may define
$$v_k(x)=u_{\alpha_k}(x_{\alpha_k}+\lambda_{\alpha_k} x).$$
It is well-known that $v_k$ converges in $C^\infty(D_R)$ to a harmonic map $v^1:\C\rightarrow N$
for any fixed $R$, where $D_R=D_R(0)=\{z: |z|<R\}\subset\mathbb{C}$ is a disk with radius $R>0$.
$v^1$ can be regarded as a harmonic map from $S^2$ into $N$.
Usually, $v^1$ is called the first bubble. For the details to get all the bubbles
we refer to the appendix of \cite{Li-Wang1}. Moreover, in \cite{Li-Wang2}
(see also \cite{Chen-Tian, Hong-Ying}) we prove the following theorem which will be used later.

\begin{thm}\label{energy.identity.min} Let $(\Sigma,g)$ be a closed Riemann surface and $N$
a compact Riemannian manifold.  Suppose that $H$ is a fixed homotopy class of maps from $\Sigma$
into $N$ and $u_\alpha$ is a minimizer of $E_{\alpha}$ in the set
$W^{1,2\alpha}(\Sigma, N)\cap H$.  Then, when $\alpha\rightarrow 1$ there
exists a subsequence $\{u_\alpha\}$ and harmonic map $u_0$ such
that $\{u_\alpha\}$ converges to $u_0$ weakly in $W^{1,2}(\Sigma, N)$ and
blows up at finitely many points $\{p_i: i=1,2,\cdots,m\}$. Moreover,
associated with each $\{p_i\}$ there exist finitely many harmonic
maps $w_{i_j}$ from $S^2$ into $N$, $j=1,2,\cdots,i_0$, such that
 $$\lim_{\alpha\rightarrow 1}E_{\alpha}(u_\alpha)
 =E(u_0)+\sum_{i=1}^m\sum_{j=1}^{i_0}E(w_{i_j}).$$
\end{thm}

\section{\bf The construction of counter-example}

\subsection{Constructing manifold $(N,h)$}
Let $h_1$ be the standard metric on
$$Y_1=\mathbb{T}^3=S^1\times S^1\times S^1=\mathbb{R}^3/2\pi\mathbb{Z}\oplus2\pi\mathbb{Z}\oplus2\pi\mathbb{Z}.$$
We denote a geodesic ball with radius $r$ in $\mathbb{T}^3$, centered at $p$, by $B_r(p)$. Fix a point $p\in Y_1$, and set
 $$X_1=\mathbb{T}^3\setminus B_r(p),$$
where $r<\frac{\pi}{4\sqrt{3}+2}$.
It is easy to see that the injective radius of $Y_1$ at $p$ is $\pi$ and $B_\pi(p)\setminus B_r(p)$ is isometric to
$$\mathbb{T}_0=\left(S^2\times(-\log\pi,\,-\log r ],\s e^{-2t}(d\mathfrak{s}^2+dt^2)\right),$$
where $g_\mathfrak{s} = d\mathfrak{s}^2$ is the standard metric over $S^2$.
It is also easy to check that $\mathbb{T}_0$ is isometric to
$$\mathbb{T}_0'=\left(S^2\times[0,\,\log\frac{\pi}{r}),\s e^{2t+2\log r}(d\mathfrak{s}^2+dt^2)\right)$$
and
$$\mathbb{T}_0''=\left(S^2\times(-\log\frac{\pi}{r},\,0],\s e^{-2t+2\log r}(d\mathfrak{s}^2+dt^2)\right).$$

Let $(X_2,h_2)=(X_1,h_1)$.
We consider the quotient space of $X_1\cup X_2$, obtained by gluing every point $x\in\partial X_1$ with the same point
$x\in \partial X_2$ together. In this way, we get a closed compact manifold $N$ and a
projection map $\phi: X_1\cup X_2\rightarrow N$.  We set $$M=\phi(\partial B_r(p)).$$

On $N\setminus M$, the metric $h_0=(\phi^{-1})^*(h_1)\cup
(\phi^{-1})^*(h_2)$ is well-defined and can be extended
to a metric $g_0$ over $N$. However, $g_0$ is not smooth and need to be modified.
Obviously, $M$ has a neighborhood which is isometric to
$$T=\left(S^2\times(-\log\frac{\pi}{r},\,\log\frac{\pi}{r}),\s e^{2|t|+2\log r}(d\mathfrak{s}^2 + dt^2)\right).$$
In fact, $T$ is obtained by gluing $\mathbb{T}_0'$ and $\mathbb{T}_0''$ along $S^2\times\{0\}$.\\

We let $\psi$ be a smooth function defined on $(-\log\frac{\pi}{r},\, \log\frac{\pi}{r})$
which satisfies the following:
\begin{enumerate}
\item[1).] $\psi=e^{2|t|+2\log r}$ when $|t|\geq \log 2$;
\item[2).]  $\psi'<0$  on $(-\log2, 0)$ and $\psi'>0$ on
$(0,\log 2)$.
\end{enumerate}
Note that 2) implies that $0$ is the only critical point of $\psi$
on $(-\log 2,\,\log 2)$. \\

We define a new metric $h$ on $N$ which is $h_0$ on $N\setminus T$, and $\psi(t)(d\mathfrak{s}^2 + dt^2)$ on $T$.
It is easy to see that $h$ is smooth on $N$. For convenience, we set
$$Q(a)=S^2\times\left(-\log \frac{a}{r},\,\log \frac{a}{r}\right)\subset T.$$
Obviously, we have $$\phi^{-1}(Q(a))\cap X_1=B_a(p)\setminus B_r(p)\subset Y_1.$$

\begin{lem}\label{in.Q(a)} Let $(N,h)$, $T$ and $Q(a)$ be defined as in the above.
Assume that $u:S^2\rightarrow (N,h)$ is a non-trivial harmonic map with $u(S^2)\subset Q(\pi)=T$. Then
$u$ is a harmonic map from $S^2$ into $M$.
\end{lem}

\proof Let $u=(v,f):S^2\rightarrow Q(\pi)$ be a harmonic
map, where $v\in C^\infty(S^2,S^2)$ and $f\in C^\infty(S^2)$.
The energy can be written as
$$E(u)=\frac{1}{2}\int_{S^2}|\nabla u|^2dV=
\frac{1}{2}\int_{S^2}(|\nabla v|^2+|\nabla f|^2)\psi(f)dV.$$
Here $dV=dV_{g_\mathfrak{s}}$ is the standard volume form of $S^2$. By a direct calculation,
it is easy to see that $u$ satisfies the following equation:
\begin{equation}\label{equation}
\left\{\begin{array}{ll}
             -\nabla(\psi(f)\nabla v)+\psi(f)|\nabla v|^2v=0,\\
            -\nabla(\psi(f)\nabla f)+\frac{1}{2}(|\nabla v|^2+|\nabla f|^2)\psi'(f)=0.
\end{array}\right.
\end{equation}
Multiplying the both sides of the second equation of the the above system by $f$
and then integrating the obtained identity over $S^2$, we get
$$\int_{S^2}\left(|\nabla f|^2\psi(f)+ \frac{1}{2}(|\nabla v|^2+|\nabla f|^2)\psi'(f)f\right)dV=0.$$
Noting that there always holds true $\psi'(f) f\geq 0$, we infer from the above identity
$$\int_{S^2}|\nabla f|^2\psi(f)d\mathfrak{s}=\frac{1}{2}\int_{S^2}(|\nabla v|^2+|\nabla f|^2)\psi'(f)fdV=0.$$
Then, this implies that $\nabla f=0$ and $f$ is a constant. Moreover, from the above identity we also have
$$|\nabla v|^2\psi' f\equiv 0.$$
Since $u$ is not trivial by the assumption, then there always exists a point $x_1\in S^2$
such that $|\nabla v|(x_1)\neq 0$. Hence we conclude that $\psi'(f)f \equiv 0$ which implies $f\equiv0$.
Then, it follows that $v$ is a harmonic map from $S^2$ into $M$.
\endproof

\begin{lem}\label{outside.Q(2r)} Let $(N, h)$ and $Q$ be the same as in Lemma 3.1. Assume that $u$ is a
harmonic map from $S^2$ into $(N, h)$ such that $u(S^2)\cap Q(2r)\neq\emptyset$ and $u(S^2)\cap\partial Q(\pi)
\neq\emptyset$. Then, we have $$E(u)\geq \pi(\pi-2r)^2.$$
\end{lem}

\proof Without loss of generality, we assume $p_1\in X_1$ s.t. $p_1\in \partial B_\pi(p)$
in $Y_1$ and $\phi(p_1)\in u(S^2)$.  First, $u$ is a branched minimal surface since $u$
is a harmonic map from $S^2$ into $N$. On the other hand, as $h$ is flat on
$\phi(B_{\pi-2r}(p_1))$, it is easy to check that $u(S^2)\cap\phi(B_{\pi-2r}(p_1))$
is a stationary varifold. Let $\mu(u(S^2)\cap B_{\pi-2r}(p_1))$ denote the area of
$u(S^2)\cap B_{\pi-2r}(p_1)$. By the monotonicity inequality for stationary varifold
(see \cite{Simon}), we have
$$\frac{\mu(u(S^2)\cap B_{\pi-2r}(p_1))}{\pi(\pi-2r)^2}\geq 1.$$
In view of the inequality and the fact $E(u)\geq \mu(u(S^2)\cap B_{\pi-2r}(p_1))$,
we derive the desired inequality $$E(u)\geq \pi(\pi-2r)^2.$$
Thus, we complete the proof.
\endproof

Since $h$ is flat on $N\setminus Q(2r)$, we have the following

\begin{lem}\label{flat}  Let $(N, h)$ and $Q$ be the same as in Lemma 3.1. Then, there is no harmonic map $u:S^2\rightarrow (N, h)$ such that
$u(S^2)\cap \overline{Q(2r)}=\emptyset$.
\end{lem}

By the definition of $\psi$, it is easy to check that $$4\pi\psi(0)\leq 16\pi r^2<\frac{1}{3}\pi(\pi-2r)^2$$ when $r$ is small enough.
Using Lemma \ref{outside.Q(2r)} and Lemma \ref{flat}, we get the following

\begin{cor} \label{E}  Let $(N, h)$ and $Q$ be the same as in Lemma 3.1. Assume that $u$ is a non-trivial harmonic map with
$E(u)<\pi(\pi-2r)^2$, then, we have $$E(u)=4m\pi\psi(0),$$ where $m$
is a positive integer.
\end{cor}
It is easy to check that there holds true $$12\pi\psi(0)<
48\pi r^2<\pi(\pi-2r)^2,$$
if $r<\frac{\pi}{4\sqrt{3}+2}$. Therefore, we know that, if
$E(u)<12\pi\psi(0)$ and $u$ is a non-trivial harmonic map, then
$E(u)=4\pi\psi(0)$ or $8\pi\psi(0)$.\\

\subsection{Homotopy class  $[u_k]$}

We have $\pi_1(Y_1)=\pi_1(\mathbb{T}^3)=\mathbb{Z}^3$.
Let $\beta\in \pi_1(Y_1)$ which represents $(1,0,0)$.
Let $x_1$, $x_2\in M$, and $\gamma_0$ be a curve in $M$ such that
$\gamma_0(0)=x_2$, and $\gamma_0(1)=x_1$. Let
$\gamma_k:[0,1]\rightarrow X$ be a curve with
$\gamma_k(0)=x_1, \gamma_k(1)=x_2$ and $[\gamma_k+\gamma_0]=k\beta$.
Let $w_0$ be a diffeomorphism from $S^2$ onto $M$ satisfying $w_0(0, 0, 1)=x_1$
and $w_0(0, 0,-1)=x_2$, where $(0,0,1)$ and $(0,0,-1)$ are the north pole
and the south pole of $S^2\subset\mathbb{R}^3$ respectively.

For the sake of convenience, we introduce the stereographic projection coordinates on $S^2$ with the south pole corresponding to $\infty$.
Thus, $w_0:S^2\rightarrow N$ can be viewed as a map from $\mathbb{C}\cup\{\infty\}$ into $N$. For simplicity,
we neglect the the stereographic projection map $\mathfrak{S}:S^2\rightarrow\mathbb{C}\cup\{\infty\}$ and
still denote $w_0\circ\mathfrak{S}^{-1}$ by $w_0$.

By the continuity of $w_0$, there exists a small $\delta_0>0$ such that $w_0(D_{\delta_0})$ is contained
in a small neighborhood of $x_1$, where $D_{\delta_0}=\{z\in\mathbb{C}:|z|<\delta_0\}$, and there exists
a large $R_0>0$ such that $w_0(\mathbb{C}\setminus D_{R_0})$ is contained in a small neighborhood of $x_2$, where $D_{R_0}=\{z\in\mathbb{C}:|z|<R_0\}$.

In order to construction a sequence of maps, we need to define the following two smooth non-negative functions $\lambda$ and $\nu$
on $[0,\,\infty)$:

(1). $\lambda(s):[0, \,\infty]\rightarrow [0,\, 1]$ with $\lambda(s)\equiv 0$ as $s\in[0, \,\delta_0]$
and $\lambda(s)\equiv 1$ as $s\in[2\delta_0,\,\infty)$.

(2). $\nu(s):[0, \, \infty)\rightarrow [0,\, 1]$ with $\nu(s)\equiv 1$ as $s\in [0,\, R_0-R^c_0]$ and $\nu(s)\equiv 0$
as $s>R_0$, where $R^c_0$ is a small positive constant number.

Now, we define a sequence of maps $u_k: S^2\rightarrow N$ by
$$u_k=\left\{\begin{array}{ll}
w_0(\lambda(|z|)z), &|z|\geq\delta_0,\\
\gamma_k\left(\dis\frac{\log |z|-\log R_0\epsilon_0}{\log\delta_0-\log
R_0\epsilon_0}\right), &R_0\epsilon_0<|z|<\delta_0,\\
w_0\left(\displaystyle\frac{z}{\nu(\frac{|z|}{\epsilon_0})\epsilon_0}\right),&|z|
\leq \epsilon_0 R_0.
\end{array}\right.$$
Here $\epsilon_0>0$ is a fixed constant number such that $R_0\epsilon_0<\delta_0$. By the arguments in \cite{Li-Wang1},
for any $i\neq j$, $u_i$ is not homotopic to $u_j$. For the sequence $\{u_i\}$ constructed in the above,
we have the following lemma:

\begin{lem}\label{inf.E.u_k} Let $u_k$ be the maps from $S^2$ into $(N, h)$ constructed in the above and $[u_k]$ denote
the class of maps in $W^{1,2}(S^2, N)\cap C(S^2, N)$, each map of which is homotopic to $u_k$. For any fixed $k$, we have
the following
$$\inf_{u\in[u_k]} E(u)=8\pi\psi(0).$$ Moreover, the $\inf E(u)$ can not
be attained by a harmonic map belonging to $[u_k]$.
\end{lem}

\proof

First of all, we prove that, for every fixed $k$, there holds true
\begin{equation}\label{leq}
\inf_{u\in [u_k]}E(u)\leq 8\pi\psi(0).
\end{equation}
Denote $z_1=(0,0,1)$ and $z_2=(0,0,-1)\in S^2$. Without loss of generality, we assume
$w_0$ is a harmonic map from $S^2$ into $M$ with $E(w_0)=4\pi\psi(0)$ with $w_0(z_1)=x_1$ and $w_0(z_2)=x_2$.
Let $\mathfrak{S}$ be the stereographic projection from
$S^2\setminus\{z_2\}$ to $\C$ and
$$\hat{u}_0(z) = w_0(\mathfrak{S}^{-1}(z)): \mathbb{C}\cup\{\infty\}\rightarrow N.$$
Choose a coordinate system $(y_1, y_2, y_3)$ in a geodesic ball $B_\rho(x_1)$ around $x_1\in N$ with $x_1=(0,0,0)$ and
$\{(y_1,y_2,0):(y_1,y_2,0)\in B_\rho(x_1)\}\subset M$. By the continuity of $w_0$, there exists a small $\delta>0$ such
that $w_0(z)\in B_\rho(x_1)$ as $|z|<\delta$. We define
$$u_0'=\eta_1 \hat{u}_0,$$
where $\eta_1$ is a smooth non-negative function which equals 1 outside $D_{2\delta}$, 0 on $D_\delta$,
and satisfies $|\nabla\eta_1|<\frac{C}{\delta}$. Here $D_{2\delta}\subset\mathbb{C}$ denotes the disk centered at the origin.
Then, we have
$$\int_{D_{2\delta}}|\nabla u_0'|^2dx^2\leq
2\int_{D_{2\delta}}(|\nabla\eta_1|^2|\hat{u}_0|^2+|\nabla \hat{u}_0|^2)dx^2
\leq C\delta.$$
Thus $u_0'$ satisfies
$$\text{dist}^M(u_0',\hat{u}_0)<C\delta,\s E(u_0')<4\pi\psi(0)+C\delta,\s
and\s u_0'(D_{\delta})=x_1.$$

Since $E$ is conformal invariant, $\hat{u}_0(\frac{1}{z})$ is also a harmonic map
from $\C\setminus\{0\}$ into $N$ with
$$E(\hat{u}_0(\frac{1}{z}),\C)=E(\hat{u}_0(z),\C).$$
Thus, $\hat{u}_0(\frac{1}{z})$ can be extended smoothly to $\{0\}$.
Choose a coordinate system $(y_1, y_2, y_3)$ in a geodesic ball $B_\rho(x_2)$ around $x_2\in N$ with $x_2=(0,0,0)$ and
$\{(y_1,y_2,0):(y_1,y_2,0)\in B_\rho(x_2)\}\subset M$. By the continuity of $w_0$, there exists a large $R>0$ such that
$\hat{u}_0(z)\in B_\rho(x_2)$ as $|z|>R$.
Then, we have
$$\hat{u}_0(\frac{1}{z})=O(z)\s\s \text{and}\s\s |\nabla\hat{u}_0(\frac{1}{z})|=O(1),\s \text{as} \s z\rightarrow 0.$$
Hence, we have
$$\hat{u}_0(z)=O(\frac{1}{z})\s \s\text{and}\s\s |z^2\nabla\hat{u}_0(z)|=O(1), \s \text{as} \s z\rightarrow\infty.$$
Let
$$u_0''(z)=\eta_2(|z|)\hat{u}_0(z),$$
where $\eta_2(|z|)$ is a smooth non-negative function which equals 0 outside $D_{R}$, 1 on $D_\frac{R}{2}$,
and satisfies $|\nabla\eta_1|<\frac{C}{R}$.
Then, we have
$$\int_{\C\setminus D_{\frac{R}{2}}}|\nabla u_0''|^2dx^2\leq
2\int_{D_{R}\setminus D_\frac{R}{2}}(|\nabla\eta_2|^2|\hat{u}_0|^2+|\nabla
\hat{u}_0|^2)dx^2
\leq \frac{C}{R}.$$
Thus
$$\text{dist}^M(u_0'',u_0)<\frac{C}{R},\s E(u_0'')<4\pi\psi(0)+
\frac{C}{R},\s
\s\text{and}\s u_0''(\C\setminus D_{R})=x_2.$$

We define
$$\phi_k=\left\{\begin{array}{ll}
u_0'(z), &|z|\geq\delta,\\
\gamma_k\left(\dis\frac{\log |z|-\log R\epsilon}{\log\delta-\log
R\epsilon}\right), &R\epsilon<|z|<\delta,\\
u_0''\left(\frac{z}{\epsilon}\right),&|z|
\leq \epsilon R.
\end{array}\right.$$
By a direct calculation, we obtain
$$\begin{array}{lll}
  \dint_{D_\delta\setminus D_{R\epsilon}}|\nabla \phi_k|^2&=&
     2\pi\dint_{R\epsilon}^\delta\left|\frac{\partial\gamma_k}{\partial r}\right|^2rdr\\
     &<&\dis\frac{c\|\dot{\gamma}_k\|_{L^\infty}^2}
     {(-\log R\epsilon+\log\delta)^2}
\dint_{R\epsilon}^\delta\frac{dr}{r}
   =\dis\frac{c\|\dot{\gamma}_k\|_{L^\infty}^2}{\log\delta-\log{R\epsilon}}.
\end{array}$$
Thus, for any $\epsilon_1>0$, we can choose suitable $\delta$, $R$
and $\epsilon$, such that
$$E(\phi_k)<8\pi\psi(0)+\epsilon_1.$$
Obviously, $\varphi_k=\phi_k(\mathfrak{S}^{-1})$ is homotopic to $u_k$, denoted by $\varphi_k\sim u_k$.
Thus, we get \eqref{leq}.\\

Next, we prove that $\inf_{u\in[u_k]} E(u)$ can not be attained by a harmonic map. Assuming it is attained by
a harmonic map $v_0$. Recall that $$8\pi \psi(0)<12\pi\psi(0)<48\pi r^2<\pi(\pi-2r)^2,$$ where $r>0$ is small enough.
By Lemma \ref{outside.Q(2r)}, $v_0(S^2)\subset Q(\pi)$. Thus $v_0$ is a harmonic map from $S^2$ into $M$.
This contradicts the fact $v_0\sim u_k$.  Hence $\inf_{u\in [u_k]} E(u)$
can not be attained by a harmonic map.

Let $u_\alpha$ be the $\alpha$-harmonic map such that, for fixed $k$,
$$E_\alpha(u_\alpha)=\inf_{u\in[u_k]\cap W^{1,2\alpha}(S^2, N)}E_\alpha(u).$$
Then, each map of $\{u_\alpha\}$ is minimizing and belongs to $[u_k]$.
We claim that $\{u_\alpha\}$ does not converge smoothly. Otherwise, the limit map is a harmonic map from $S^2$ into $N$,
which is homotopic to $u_k$. This contradicts the above fact $\inf_{u\in[u_k]} E(u)$ can not be attained by a harmonic map.
Hence, the bubbles must appear in the convergence of $u_\alpha$. If we denote the weak limit of $\{u_\alpha\}$ as $u_0$  and
the bubbles as $v^1$, $\cdots$, $v^m$, then, by Theorem \ref{energy.identity.min}, we have
$$\inf_{u\in [u_k]} E(u)=\lim_{\alpha\rightarrow 1}E_{\alpha}(u_\alpha)=
E(u_0)+\sum_{i=1}^m E(v^i).$$
Since $E(u_0)$ and $E(v^i)$ are smaller than $\pi(\pi-2r)^2$, then, we know that
$E(u_0)+\sum_{i=1}^m E(v^i)$ equals only $8\pi\psi(0)$ or $4\pi\psi(0)$.\\

Next, we will show that the following identity does not hold true
$$E(u_0)+\sum_{i=1}^mE(v^i) = 4\pi\psi(0).$$
If this is true, then $u_0$ is trivial and $u_\alpha$ has only one bubble $v^1$.
To derive a contradiction, we only need to prove $u_\alpha\sim v^1$.

Let $x_0\in S^2$ be a blowup point. Take an isothermal coordinate system around $x_0$
with $x_0=(0,0)$ on $S^2=\C\cup\{\infty\}$. Let $v^1$ is the limit map of
$u_\alpha(z_\alpha+\lambda^1_\alpha z)$, where $z_\alpha\rightarrow 0$,
$\lambda^1_\alpha\rightarrow 0$.
Then, we know that $$v^1_\alpha(z)=u_\alpha(z_\alpha+\lambda^1_\alpha z)$$ converges smoothly to $v^1$
on any $D_R=D_R(0)\subset\mathbb{C}$. Moreover, $u_\alpha$ converges smoothly in $\C\cup\{\infty\}\setminus
D_\frac{1}{R}$ to a point $y_0\in N$. For us to prove $u_\alpha\sim v^1$,
it is enough to check that for any $\epsilon>0$, there exists an $R>0$, such that
$$\sup_{t\in [R\lambda^1_\alpha,\,\frac{1}{R}]}
\text{osc}_{\partial D_{t}(z_\alpha)}
u_\alpha<\epsilon.$$
Indeed, if this is not true,  then, there exists a sequence of $\lambda^2_\alpha$
with $\lambda^2_\alpha\rightarrow 0$ and $\lambda^2_\alpha/\lambda^1_\alpha\rightarrow+\infty$,
such that
$$\text{osc}_{\partial D_{\lambda^2_\alpha}(z_\alpha)}
u_\alpha
\rightarrow \epsilon_1\neq 0.$$
Let $$v^2_\alpha(z)=u_\alpha(z_\alpha+\lambda^2_\alpha z).$$ If the sequence $\{v^2_\alpha\}$
has blowup points, then, at each blowup point there exists at least a bubble of
$\{v^2_\alpha\}$, which is also a bubble of $u_\alpha$ and is different from the previous bubble $v^1$.
However, it is impossible since there only exists one bubble for $\{u_\alpha\}$. Hence, we infer that,
as $\alpha\rightarrow 1$, $\{v^2_\alpha\}$
converges smoothly on $D_{R'}\setminus D_{\frac{1}{R'}}\subset\mathbb{C}$ for any $R'$. It follows that
$$\text{osc}_{\partial D_1}v^2_\alpha
\rightarrow \epsilon_1\neq 0.$$
This means that the limit map of $\{v^2_\alpha\}$ is not trivial and the limit map is also a bubble
of $\{u_\alpha\}$ which is different from $v^1$. This is a contradiction. Thus, we conclude
$$\inf_{u\in [u_k]} E(u)=E(u_0)+\sum_{i=1}^m E(v^i)=8\pi\psi(0).$$
Thus, we complete the proof of the lemma.
\endproof

By the Sobolev embedding theorem we know that, as $\alpha>1$, $$W^{1,2\alpha}(S^2, N)\subset C(S^2, N).$$
For simplicity, let $[u_k]^\alpha$ denote the the class of maps belonging to $W^{1,2\alpha}(S^2, N)$,
each map of which is homotopic to $u_k$. In fact, it is easy to see that
 $$[u_k]^\alpha=[u_k]\cap W^{1,2\alpha}(S^2, N).$$
From now on, we always let $u_{\alpha,k}$ denote the smooth map which attains $\inf_{u\in [u_k]^\alpha} E_\alpha(u)$, i.e.,
$$E_\alpha(u_{\alpha,k})=\inf_{u\in[u_k]^\alpha} E_\alpha(u).$$

\begin{lem} For any $\lambda_0>8\pi\psi(0)$, there always exists a sequence $\{\alpha_k\}$ with $\alpha_k\rightarrow 1$ and a sequence
$\{i_k\}$ such that $E_{\alpha_k}(u_{\alpha_k,i_k})= \lambda_0$ for every $k$.
\end{lem}

\proof For $\alpha\in [1,\, \alpha_0)$ where $\alpha_0-1>0$ is small enough, we define the following function
$$\varphi_k(\alpha)=\inf_{u\in [u_k]^\alpha} E_\alpha(u).$$
Firstly, we need to show that, for any fixed $\alpha\in(1,\alpha_0)$,
\begin{equation}\label{infinity}
\lim_{k\rightarrow+\infty}
\varphi_k(\alpha)=+\infty.
\end{equation}
If this is false, then, there exists a constant $C$ such that $\varphi_k(\alpha)\leq C$ as $k$ is large enough.
We note that there holds true, for any small $\delta$ and $x\in S^2$,
\begin{equation}\label{33}
E(u_{\alpha,k},B_\delta(x))=\frac{1}{2}
\int_{B_\delta(x)}|\nabla u_{\alpha,k}|^2\leq
\frac{1}{2}\left(\int_{B_\delta(x)}|\nabla u_{\alpha,k}|^{2\alpha}\right)^\frac{1}{\alpha}
|B_\delta(x)|^\frac{\alpha-1}{\alpha}.
\end{equation}
Hence, we can pick a fixed $\delta$, which is small enough, such that
$$E(u_{\alpha,k},B_\delta(x))<\epsilon_0.$$
Thus, by Theorem \ref{epsilon}, there exists a subsequence of $u_{\alpha,k}$ which converges smoothly to a smooth map $u_0$
as $k$ tends to $\infty$. Hence, we know that $u_{\alpha,k}$ are homotopic to $u_0$ for any $k$.  This contradicts
the fact that $u_{\alpha,i}$ is not homotopic to $u_{\alpha,j}$ as $i \neq j$.

Next, we want to prove $\varphi_k$ is continuous on $[1,\alpha_0)$.
Using \eqref{33} again, we can prove that, for a fixed small $\epsilon>0$,
there holds true $$\|\nabla u_{\alpha,k}\|_{C^0(S^2)}<\Lambda(\epsilon)$$ for any
$\alpha\in (1+\epsilon,\alpha_0)$.
For any $\alpha$, $\alpha'\in(1+\epsilon,\alpha_0)$, we have
$$\varphi_k(\alpha)
\geq\frac{1}{2}(1+C_1^2)^{\alpha-\alpha'}
\int_{S^2}(1+|\nabla u_{\alpha,k}|^2)^{\alpha'}-\frac{1}{2}
\geq (1+C_1^2)^{\alpha-\alpha'}\varphi_k(\alpha')+
\frac{1}{2}(1+C_1^2)^{\alpha-\alpha'}-\frac{1}{2},$$
where
$$C_1=\left\{\begin{array}{ll}
             0,&\text{when}\s \alpha>\alpha',\\
            \Lambda(\epsilon),&\text{when}\s \alpha<\alpha'.
\end{array}\right.$$
It follows that
$$\varliminf_{\alpha\rightarrow \alpha'}\varphi_k(\alpha)\geq\varphi_k(\alpha').$$
On the other hand, we also have
$$\varphi_k(\alpha')\geq
\frac{1}{2}(1+C_2^2)^{\alpha'-\alpha}
\int_{S^2}(1+|\nabla u_{\alpha',\, k}|^2)^{\alpha}-\frac{1}{2},$$
where
$$C_2=\left\{\begin{array}{ll}
             0, &\text{when}\s \alpha'>\alpha,\\
            \|\nabla u_{\alpha',\, k}\|_{L^\infty}, &\text{when}\s\alpha'<\alpha.
\end{array}\right.$$
It follows that
$$\varphi_k(\alpha')\geq
(1+C_2^2)^{\alpha'-\alpha}
\varphi_k(\alpha)
+\frac{1}{2}(1+C_2^2)^{\alpha'-\alpha}
-\frac{1}{2}.$$
Hence, we have
$$\varlimsup_{\alpha\rightarrow\alpha'}\varphi_k(\alpha)\leq \varphi_k(\alpha').$$
Therefore, we have
$$\lim_{\alpha\rightarrow\alpha'}\varphi_k(\alpha)=\varphi_k(\alpha').$$
Hence, we show the continuity of $\varphi_k(\alpha)$ on $(1,\,\alpha_0)$.

Next, we want to prove that $\varphi_k(\alpha)$ is left continuous at 1. Equivalently, we need to show
\begin{equation}\label{limit}
\lim_{\alpha\searrow 1}\varphi_k(\alpha)=\varphi_k(1).
\end{equation}
Obviously, we have that, for any fixed $u\in W^{1,2}(S^2, N)$ and $\alpha_1>\alpha_2>1$,
$$E_{\alpha_1}(u)\geq E_{\alpha_2}(u)\geq E(u).$$
It follows that
$$\varphi_k(\alpha_1)\geq \varphi_k(\alpha_2)\geq \varphi_k(1).$$
Hence, $\lim_{\alpha\searrow 1}\varphi_k$ exists and
$$\lim_{\alpha\searrow 1}\varphi_k(\alpha)\geq \varphi_k(1).$$

On the other hand, note $u_k$ is a smooth map. Then, for any $\epsilon>0$, there exists
a smooth map $u'_k\in C^\infty(S^2,N)$ which is homotopic to $u_k$, i.e., $u'_k\sim u_k$,
and satisfies
$$E(u'_k)\leq \varphi_k(1)+\epsilon.$$
Since
$$\lim_{\alpha\searrow 1}E_{\alpha}(u_k')=E(u'_k)\s\s \text{and}\s \s
\varphi_k(\alpha)\leq E_{\alpha}(u_k'),$$
we have
$$\lim_{\alpha\searrow 1}\varphi_k(\alpha)\leq \varphi_k(1)
+\epsilon,$$
which implies \eqref{limit}. Hence, we show that $\varphi_k(\alpha)$
is continuous on $[1,\,\alpha_0)$ for any fixed $k$.\\

By \eqref{infinity}, for any given sequence $\{\alpha_k'\}$ with $\alpha_k'\rightarrow 1$,
there exists always a sequence $\{i_k\}$, such that $E_{\alpha'_k}(u_{\alpha'_{k},i_k})>\lambda_0$, i.e.,
$\varphi_{i_k}(\alpha'_k)>\lambda_0$. Lemma 3.5 tells us that $\varphi_{i_k}(1)=8\pi\psi(0)$ for any $i_k$.
By the assumption $\lambda_0>8\pi\psi(0)$ we have
$$\varphi_{i_k}(\alpha'_k)>\lambda_0>\varphi_{i_k}(1).$$
Since $\varphi_k(\alpha)$ is continuous on $[1,\,\alpha_0)$,
we conclude that, for any fixed $i_k$, there always exists $\alpha_k\in (1,\,\alpha'_k)$ such that
$$\varphi_{i_k}(\alpha_k)=E_{\alpha_k}(u_{\alpha_k,i_k})=\lambda_0.$$ Thus, the proof is finished.
\endproof

\subsection{The counter-example}
By Lemma 3.6, for given $\tau\in (8\pi\psi(0),12\pi\psi(0))$ there exist a sequence
$\{\alpha_k:\alpha_k>1, k\in\mathbb{N}\}$ with $\alpha_k\rightarrow 1$ and a sequence
of minimizing $\alpha_k$-harmonic maps $v_k\in W^{1, 2\alpha_k}(S^2, N)$ with $v_k\sim u_{i_k}$
such that
$$\tau=E_{\alpha_k}(v_k)=\inf_{u\in [u_{i_k}]^{\alpha_k}}E_{\alpha_k}(u),\s\s \forall k\in \mathbb{N}.$$
Since $v_i$ and $v_j$ are not in the same homotopy class for any $i\neq j$, we known that $v_k$
must blows up as $k\rightarrow+\infty$. Let $v^0$ be the weak limit of $\{v_k\}$ in $W^{1,2}(S^2, N)$,
and $v^1$, $\cdots$, $v^m$ are all the bubbles produced in the convergence of $\{v_k\}$. Since
$E(v^i)<12\pi\psi(0)$, then, it follows Corollary 3.4 that $E(v^i)=4\pi\psi(0)$ or $8\pi\psi(0)$.
Hence, $$\frac{1}{4\pi\psi(0)}(E(v^0)+\sum_{i=1}^m  E(v^i))$$ is always an integer. However, certainly
$\frac{\tau}{4\pi\psi(0)}$ is not an integer by the previous assumption. So, the energy identity
is not true for the sequence $\{v_k\}$, i.e.,
$$\lim_{k\rightarrow\infty}E_{\alpha_k}(v_k)\neq E(v^0)+\sum_{i=1}^m E(v^i).$$

\begin{rem}
By an argument in \cite{Li-Wang1}, we also have
$$\lim_{k\rightarrow\infty}E(v_k)\neq E(v^0)+ \sum_{i=1}^mE(v^i).$$
\end{rem}

\section{\bf An example of manifold with non-discrete set of energy for harmonic 2-spheres}

In this section, we will construct a Riemannian manifold $(N, h)$ for which $\mathcal{E}(N,h)$ is not discrete.
In other words, $\mathcal{E}(N,h)$ admits limit points.

Let $\psi(t)$ be a smooth positive function defined on $(-1,\,1)$ satisfying
$$\psi(t)=e^{-\frac{1}{t^2}}\sin\frac{1}{t}+1,\quad\quad t\in (-\frac{1}{2}, \, \frac{1}{2}).$$
It is easy to check that the critical point of $\psi(t)$ satisfies the equation
$$\tan \frac{1}{t}=\frac{t}{2}.$$
Thus, we can find $t_k\rightarrow 0$, such that $\psi'(t_k)=0$, $\psi(t_k)\neq 1$
and $\psi(t_k)\rightarrow 1$.

Let
$$h=\psi(t)(d\mathfrak{s}^2+dt^2),$$
which is a metric over $S^2\times(-1,1)$.  Let
$v$ be the identity map from $S^2$ to $S^2$ and
$$u_k=(v,t_k): S^2\rightarrow (N, h)\equiv (S^2\times(-1,1), h).$$
By \eqref{equation}, it is easy to see that $u_k$ is a harmonic map from $S^2$ into $(S^2\times(-1,1),h)$
with
$$E(u_k)=4\pi\psi(t_k).$$
Thus, $4\pi$ is not a discrete number in $\mathcal{E}(S^2\times(-1,1),h)$.\\

\noindent{\bf Acknowledgement:} {\small The authors are grateful to thank Professor Weiyue Ding for his encouragement
and beneficial suggestions. The first author is supported by NSFC, Grant No. 11131007; the second author is supported by NSFC, Grant No. 10990013.}

\end{document}